\documentclass[11pt]{article}

\usepackage{url}
\usepackage{booktabs}
\usepackage{tabularx}
\usepackage{float}             
\usepackage{graphicx}
\usepackage{amsfonts}
\usepackage{amsmath}
\usepackage{amsthm}
\usepackage{amssymb}
\usepackage{bm}                     
\usepackage{color}
\usepackage{indentfirst}
\usepackage{multirow}
\usepackage{longtable}
\usepackage{bbm}
\usepackage[all]{xy}
\usepackage{pdflscape}
\usepackage{authblk}
\usepackage{ulem}

\usepackage{algorithm}
\usepackage{algorithmic}
\usepackage{amsfonts}
\usepackage{amsmath}
\usepackage{amsopn}
\usepackage{amssymb}
\usepackage{bm}
\usepackage{color}
\usepackage{enumerate}
\usepackage{epstopdf}
\usepackage{graphicx}
\usepackage{lipsum}
\usepackage{mathrsfs}
\usepackage{url}

\usepackage[paper=a4paper,dvips,top=3cm,left=2.4cm,right=2.4cm,
    foot=1cm,bottom=3cm]{geometry}
 
\usepackage{array}
\newcommand{\PreserveBackslash}[1]{\let\temp=\\#1\let\\=\temp}
\newcolumntype{C}[1]{>{\PreserveBackslash\centering}p{#1}}
\newcolumntype{R}[1]{>{\PreserveBackslash\raggedleft}p{#1}}
\newcolumntype{L}[1]{>{\PreserveBackslash\raggedright}p{#1}}

\allowdisplaybreaks[4]

\newtheorem{Theorem}{Theorem}[section]

\newcommand{\E}{\acute{E}}

\begin{document}

\title{A Convergent Semi-Proximal Alternating Direction Method of Multipliers 
    for Recovering Internet Traffics from Link 
    Measurements}
\author[1,2]{Zhenyu Ming}
\author[1]{Liping Zhang \thanks{Corresponding author: 
lipingzhang@mail.tsinghua.edu.cn}}
\author[1]{Hao Wu}
\author[2]{Yanwei Xu}
\author[2]{Mayank Bakshi}
\author[2]{Bo Bai}
\author[2]{Gong Zhang}
\affil[1]{Department of Mathematical Sciences, Tsinghua University, Beijing 
100084, China}
\affil[2]{Theory Lab, Central Research Institute, 2012 Labs, Huawei 
Technologies Co., Ltd, Hong Kong}

\date{}
\maketitle

\begin{abstract}
    It is challenging to recover the large-scale internet traffic data purely 
    from the link measurements. With the rapid growth of the problem scale, it 
    will be extremely difficult to sustain the recovery accuracy and the 
    computational cost. In this work, we propose a new Sparsity Low-Rank 
    Recovery (SLRR) model and its Schur Complement Based semi-proximal 
    Alternating Direction Method of Multipliers (SCB-spADMM) solver. Our 
    approach distinguishes itself mainly for the following two aspects. First, 
    we fully exploit the spatial low-rank property and the sparsity of traffic 
    data, which are barely considered in the literature. Our model can be 
    divided into a series of subproblems, which only relate to the traffics in 
    a certain individual time interval. Thus, the model scale is significantly 
    reduced. Second, we establish a globally convergent ADMM-type algorithm 
    inspired by [Li et al., Math. Program., 155(2016)] to solve the SLRR model. 
    In each iteration, all the intermediate variables' optimums  can be 
    calculated analytically, which makes the algorithm fast and accurate. 
    Besides, due to the separability of the SLRR model, it is possible to 
    design a parallel algorithm to further reduce computational time. According 
    to the numerical results on the classic datasets Abilene and GEANT, our 
    method achieves the best accuracy with a low computational cost. Moreover, 
    in our newly released large-scale Huawei Origin-Destination (HOD) network 
    traffics, our method perfectly reaches the seconds-level feedback, which 
    meets the essential requirement for practical scenarios.
    
    {\bf keywords:} Large-scale network recovery,  HOD dataset, Spatial 
low-rankness, Nuclear norm minimization,  semi-proximal ADMM
    \end{abstract}

\section{Introduction}

The increasing demand of various services and applications running on the 
internet has led to an exponential growth of internet traffic in the last two 
decades. 
Such growth has, in turn, resulted in significant challenges for network 
resource orchestration, capacity planning, service provisioning, and traffic 
engineering.
One of the most important input factors in these tasks is traffic data 
\cite{cah1998,for2002,lak2004,ro2012,TR13,xie2019}. It comprises the volumes of 
traffics in bytes, packets, or flows during a specified period of time between 
Origin and Destination (OD) pairs which can be edge routers in the WAN network 
or the ToR switches in the data center settings in the network 
\cite{tmverdi96}. In practice, traffic data 
is often arranged in a matrix or a tensor, called the traffic 
matrix \cite{TMLink03,zhang2005} or the traffic tensor 
\cite{xie2018,xie2015,xie2019,zh2015}, respectively.

In general, there are two ways of recovering the whole OD traffics. One is to 
use completion techniques by directly monitoring part of the OD traffics. 
Extensive works have been established on this field in recent years 
\cite{ca2010,ch2012,ji2020,recht2010,ro2012,xie2017,2019Active,xie2018,xie2015,xie2019,zh2015}.
Another is to infer the OD traffic data  from pure  link-load measurements, or 
from the combination of link-loads and partial OD traffics. This is the 
Network Tomography Problem (NTP)   
\cite{TMLink04,recht2010,ro2012,tmverdi96,tmcao00,TMLink03}. In this paper, we 
aim to recover OD traffics  purely by the link-load measurements. 
It should be noted that pure completion approaches 
\cite{ji2020,la2015,xie2017,2019Active,xie2019Ele,xie2018,xie2015,xie2019} are, 
to a great extent, different from ours since such models do not refer to any 
information of link-loads. Moreover, from a practical viewpoint, the cost of 
direct measurement of OD traffics is much higher than that of link traffics. 
Thus, NTP is a more significant issue to be addressed in real-world 
network engineering.

In recent years, only a few studies have been devoted to NTP, especially for 
large-scale cases. 
The Gravity model 
\cite{zhang2005} is a rough rank-1 approximation to traffic matrix, 
using regularization based on entropy penalization. Such method is 
applicable to the cases where OD traffics are stochastically independent.
Baseline Approximation (BA) \cite{ro2012} is a rank-2 interpolation to the 
traffic matrix, with the idea of centering the traffic data. The above two 
methods often result in a rough estimation of the traffic matrix, but such 
estimation might be viewed as a warm-start of other methods. Tomo-Gravity 
\cite{zhang2005} model is another regularization model based on 
Kullback-Leibler divergence between the gravity model and the direct 
measurements. By employing the directly measured traffics, the accuracy of the 
Tomo-Gravity model has a significant improvement compared to that of the 
primitive Gravity model.
Recht et al. \cite{recht2010} studied the  minimum-rank solutions of linear 
matrix equations. The SDPLR model and the method of multipliers proposed in 
\cite{recht2010} is also 
applicable to NTP.  This research shows 
that the minimum rank solution over the given affine space can be efficiently 
recovered using nuclear norm relaxation technique if a certain restricted 
isometry property holds for the linear mapping involved in the constraints.
Sparsity Regularized Matrix Factorization (SRMF) proposed by Roughan et 
al. \cite{ro2012} is a classic approach and regarded as the benchmark in this 
field.  SRMF exploits both the global low-rank property and the spatio-temporal 
structure of traffic matrix.
They further combine their SRMF with previous 
Tomo-Gravity, and name this 
hybrid method as Tomo-SRMF. Tomo-SRMF achieved the best 
estimation of traffics compared to the others at the time. However, it does not 
consider the spatial low-rank property or the high sparsity, which will be 
introduced soon, of real-world traffic scenarios. Therefore, the recovery 
accuracy of Tomo-SRMF is still not overwhelming. 

Here we consider the spatial low-rank property. In this case, the traffic 
matrix is a square matrix whose rows and columns correspond to all of the 
origins and destinations in numerical order. Thus, it only restores the 
traffics of an individual time interval. It is called the spatial matrix. Its 
low-rank property can be easily exploited. Another matrix model for describing 
traffic is to use both spatial and temporal information. The rows and the 
columns of this matrix correspond to different OD pairs and time intervals. 
Correspondingly, this matrix has Spatio-temporal low-rank property in practice. 
Finally, when traffic data is arranged into tensor, there is also a 
corresponding low-rank model, which depends on the definition of tensor ranks, 
e.g., CP-rank \cite{Ki2000}, Tucker-rank \cite{Tu1966,Ku1989}, Tubal-rank 
\cite{Ki2011}, TT-rank \cite{Os11}. In general, this type of model is more 
complicated and not easy to implement.

Sparsity is another important feature in large-scale practical problems. This 
property arises because the number of destinations of traffics initiated from a 
specific origin is always quite smaller than the whole number of nodes 
\cite{TR13}. For instance, in a Data Center Interconnection network of a large 
online social network or resource distribution website, most of the traffics 
are the backup traffics from the Data Centers providing services for customers 
to the backup Data Centers \cite{B42018}. Besides, in the Metropolitan Area 
Networks, the traffics are always between the large Data Centers providing 
internet services and the small Data Centers of Internet Service Providers 
directly connecting residential or office areas \cite{DCTraffic2010}. On the 
other hand, only the elephant flows, especially the top-k elephant flows, are 
interested by the network administrators or network orchestrating applications 
such as congestion control, anomaly detection, and traffic engineering 
\cite{xie2019Ele}. By employing some efficient network measurement tools such 
as the Bloom Filters \cite{BF2005} or analyzing the routing tables, we can 
determine the exact positions of the zero OD pairs.

\begin{figure}[ht]
    \centerline{\includegraphics[width=13.5cm,height=4.2cm]{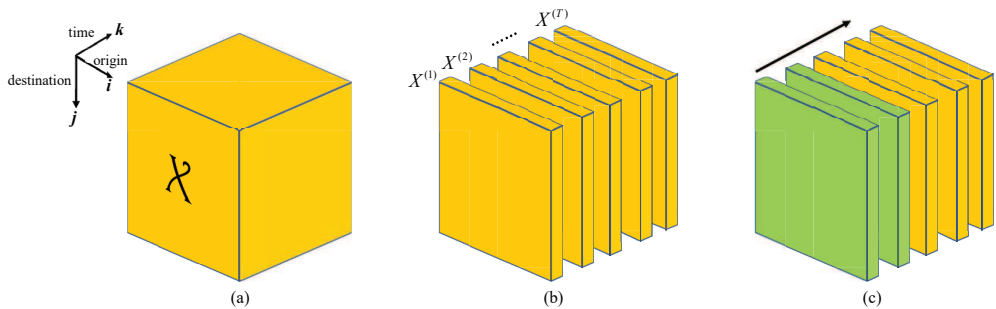}}
    \caption{Third-order traffic tensor model. (a) Third-order traffic tensor 
        $\mathcal{X}$, with the $(i,j,k)$-th element $x_{ijk}$ being the 
        traffic 
        flow transmitted from the $i$-th origin to $j$-th destination in the 
        $k$-th 
        time interval. (b) Different  time intervals of $\mathcal{X}$, i.e., 
        $X^{(1)},\cdots,X^{(T)}$. (c) The traffic tensor $\mathcal{X}$ is 
        recovered 
        as chronological order. The traffics in green have been recovered, 
        while 
        the traffics in yellow have not.}\label{tensor}
\end{figure}

In this paper, we arrange the traffic data as a third-order traffic tensor   
$\mathcal{X}(=x_{ijk})$, with each component $x_{ijk}$ being the traffic flow 
transmitting from the $i$-th node to the $j$-th node over the $k$-th time 
interval (see Figure \ref{tensor} (a)). 
Consider the network consists of S nodes and the traffic information on T 
discrete time intervals, the dimension of traffic tensor $\mathcal{X}$ is 
$S\times S\times T$, and the dimension of its frontal slices, i.e. 
$X^{(1)},\cdots,X^{(T)}$ (see Figure \ref{tensor} (b)), are all $S\times 
S$. If we further assume that the network has 
$M$ links and $N$ OD pairs, then $M,\ S$ and $N$ satisfy the following 
relations:
$$N^{1/2}=S\le M+1<N.$$ 
The character of a large-scale network system is that the size of $M$ and $N$ 
could be very large. In addition, the network topology would be rather sparse 
($S\approx M\ll N$). In the HOD network, there 
are 577 links, 243 nodes and 59049 OD pairs. More than 95\% of OD pairs do not 
transmit any traffic.

In this work, we propose a fast and accurate recovery method for NTP. First, we 
develop a convex Sparsity Low-Rank Recovery (SLRR) model, which is the first to 
consider both spatial low-rank property and the sparsity of traffic data. 
Moreover, SLRR 
can be divided into a group of subproblems, each of which only related to the 
traffics in a single time interval. Thus, the problem scale can be sharply 
reduced, and we can recover the traffics as chronological order (see 
Figure \ref{tensor} (c)). 
Then, we design a Schur Complement Based semi-proximal Alternating Direction 
Method of Multipliers (SCB-spADMM) to solve SLRR. For each block variable, we 
derive the explicit form of its optimum in the updating schemes. Thus, the 
algorithm iteration is rather fast and accurate. Based on Li, Sun and Toh's 
research \cite{li2014}, the global convergence of SCB-spADMM is established. In 
addition, with the separable structure of SLRR, people can easily design the 
parallel algorithm. Finally, we utilize cross-validation, a commonly used 
strategy for parameter setting in machine learning, to previously pick out a 
group of appropriate parameters involved in our model and algorithm. Such a 
technique is data-independent, and thus it strengthens the robustness of our 
approach.

To test the accuracy and efficiency of numerical methods, we need 
practical network datasets. In this field, Abilene 
dataset\footnote{\url{https://www.cs.bu.edu/fac/crovella/abilene-distro.tar}} 
and 
GEANT dataset\footnote{\url{https://totem.info.ucl.ac.be/dataset.html}} are 
frequently used. 
Abilene/G$\rm\E$ANT dataset consists of 41/74 links and 11/23 nodes, collected 
using 5-minute/15-minute time intervals from Dec. 8, 2003 to Dec. 28, 2003/Jan. 
1, 2005 to Apr. 29, 2005.
However, the scale of these datasets is too small compared to the real network 
problems. Therefore, we released a new dataset -- HOD 
dataset\footnote{\url{https://gitee.com/network-user123/hod-network_traffic_data}},
which is extracted and modified from real data. The traffics of HOD network are 
collected in 671 time-intervals, each of which lasts for 15 minutes. There are 
243 nodes and 577 links, and 95\% of OD pairs never transmit traffic. Thus HOD 
traffic data is highly sparse. Since monitoring the OD traffics is 
prohibitively 
expensive, our HOD dataset does not contain such information. This means that 
HOD dataset is mainly applicable to NTP that this work cares about. We will see 
that the method proposed in this work is very effective for HOD dataset in 
Section \ref{hod}. 
Furthermore, we also hope that the newly proposed HOD dataset can be used as 
the benchmark of traffic engineering methods in the future.

The outline of this paper is as follows. 
Our recovery model SLRR and 
SCB-spADMM algorithm are discussed in Section \ref{model} and Section 
\ref{algorithm}, respectively. In Section \ref{ns}, we conduct numerical 
simulations and evaluate the performance of different approaches. In Section 
\ref{hod}, we introduce the newly published HOD dataset and present the 
corresponding numerical simulations. In Section \ref{conc}, the conclusion of 
our work 
is drawn.


\section{Recovery Model}\label{model}
In this section, we first establish a regularized nuclear norm minimization 
with linear constraints, called Sparsity Low-Rank Recovery (SLRR) model, to 
infer the 
traffic data from link measurements. After analyzing the structure property of 
traffic data, we handle the traffic tensor $\mathcal{X}$ according to its 
frontal slices	$X^{(1)},\cdots,X^{(T)}$. This improvement not only notably 
reduces the size of the problem but also makes it possible to design parallel 
algorithm.
Then, we derive the dual model of SLRR, which can be easily solved by 
SCB-spADMM. 
Finally, the optimality conditions of the primal and dual problems are 
presented for setting the stopping criterion of SCB-spADMM.

\subsection{Sparsity Low-Rank Recovery Model}

\begin{figure*}[ht]
    \centerline{\includegraphics[width=13.5cm,height=3cm]{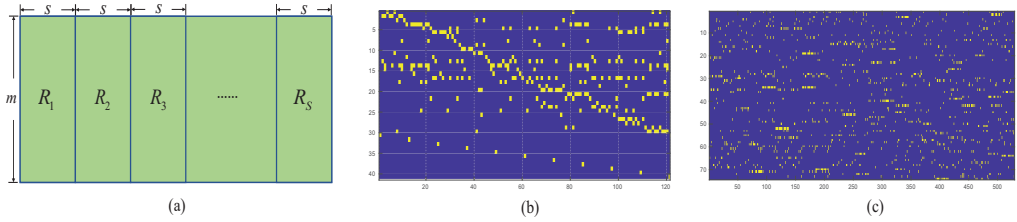}}
    \caption{Routing matrix and its visualization. (a) Routing matrix $R$ 
        and its blocks $R_1,\cdots,R_S$. The dimensions of $R$ and each block 
        $R_i(i\in\{1,\cdots,S\})$ are $m\times S^2$ and $M\times S$. (b) The 
        routing matrix of Abilene dataset. (c) The routing matrix of 
        G$\rm\E$ANT dataset. In (b) and (c), the blue elements are 0 while the 
        yellow element are 1.}\label{routing}
\end{figure*}

Intuitively, the optimization model to recover the traffics in the $k$-th 
($k=1,\cdots,T$) time 
interval $X^{(k)}$ with the spatial low-rank property can be established by 
minimizing the rank function:
\begin{equation}\label{model1}\begin{aligned}
\min\limits_{X^{(k)}} \quad&\text{rank}(X^{(k)}) \\
\text { s.t.} \quad& \mathcal{A}_k(X^{(k)})=B_k,
\end{aligned}\end{equation}
where $\mathcal{A}_k$ is a linear operator and $B_k$ contains measurements and 
other given constrains. Since rank function is nonconvex and nonsmooth, it is 
often relaxed by nuclear norm $\|\cdot\|_{*}$. It stands for the sum of all the 
singular values of a matrix, which is proved to be the tightest convex 
relaxation of 
rank function on the unit sphere $\mathcal{B}:=\{X:\ \|X\|_2\le1\}$  
\cite{fazel2002}. The linear constraint $\mathcal{A}_k(X^{(k)})=B_k$ 
contains three types of information:
\begin{enumerate}
    \item $X^{(k)}$ is nonnegative.
    \item The sparsity constraint of the complete traffic tensor takes 
    form as  
    $\mathcal{P}_{\Omega}(\mathcal{X})=0$, where $\Omega\subset\{(i,j,k):\ 
    1\le 
    i,j\le S,\ 
    1\le k\le T\}$ is the index set indicating the location of these 
    zero-traffics 
    and $\mathcal{P}_{\Omega}(\cdot)$ is the projection on $\Omega$, i.e., 
    $\mathcal{P}_{\Omega}(x_{ijk})=x_{ijk}$ if $(i,j,k)\in\Omega$, otherwise, 
    $\mathcal{P}_{\Omega}(x_{ijk})=0$.
    Denote $\Omega^{(k)}:=\{(i,j):\ (i,j,k)\in\Omega\}$.
    Thus, in model \eqref{model1}, it follows that
    $P_{\Omega^{(k)}}(X^{(k)})=0$. 
    \item By link-load measurements, we have 
    \begin{equation}\label{Linkload}
    L^{(k)}=R\cdot\text{vec}(X^{(k)}),
    \end{equation}
    where $R\in\mathbb{R}^{M\times N}$ is the 0-1 routing matrix, indicating 
    which link associates to which OD pairs (see Figure \ref{routing}),
    $L^{(k)}\in\mathbb{R}^{M}$ is the $k$-th column of link measurements matrix 
    $L$ and $\text{vec}(X^{(k)})$ is vectorization of the spatial matrix 
    $X^{(k)}$. Formulation \eqref{Linkload} can be easily rewritten as 
    $$R\cdot\text{vec}(X^{(k)})=\sum\limits_{j=1}^{S}R_jX^{(k)}e_j,$$ where 
    $R=[R_1,R_2,\cdots,R_S]$ is uniformly divided into $S$ parts along its 
    columns and $e_j$ is the zero-vector except the $j$-th element being 1.
\end{enumerate}  
Therefore, we modify model \eqref{model1} to (for simplicity but with a little 
ambiguity, we omit the index `(k)'):
\begin{equation*}\label{model2}\begin{aligned}
\min\limits_{X} \quad&\|X\|_{*} \\
\text { s.t.} \quad& \sum\limits_{j=1}^{S}R_jXe_j=L,\quad 
\mathcal{P}_{\Omega}(X)=0,\quad X\ge0.
\end{aligned}\end{equation*}

In internet traffic data, the variation of the traffics from the start to the 
end is often considerably large. Nevertheless, the traffics change slowly 
between adjacent time intervals. This continuity property of traffics has been 
studied in many models \cite{ro2012,xie2018,zh2015}. On the other hand, the 
traffics is also has periodic \cite{xie2018}. This situation of similar traffic 
in the same time period at different times can be understood as Internet users 
having similar behaviors in the same time period every day. By adding two 
regularized terms in model (2), we have obtained a more realistic model:
\begin{equation}\label{model3}\begin{aligned}
\min\limits_{X} \quad&\|X\|_{*}+\rho_1\|X-\bar{X}\|_F^2+\rho_2\|X-\hat{X}\|_F^2 
\\
\text { s.t.} \quad& \sum\limits_{j=1}^{S}R_jXe_j=L,\quad 
\mathcal{P}_{\Omega}(X)=0,\quad X\ge0.
\end{aligned}\end{equation}
Here, $\bar{X}$ and $\hat{X}$ represent the traffics one time interval ago, and 
one week ago, respectively. The regularization parameters $\rho_i > 0\ (i=1,2)$ 
allow a tunable tradeoff between the spatial low-rank property, the continuity, 
and the periodicity of traffics.

Equation \eqref{model3} leads to the Sparsity Low-Rank Recovery (SLRR) model. 
It allows us to cope with the traffics in each individual time interval as the 
chronological order $X^{(1)}\to X^{(2)}\to\cdots$, instead of the aggregate 
data. 
As shown in Figure \ref{tensor}(c), the green frontal slices represent the 
recovered traffic data, while the yellow ones are not. This character sharply 
reduces the variable size of each independent subproblem, and it correctly 
describes the basic requirement that the update traffics needs to be reported 
immediately each time when the lastest time slice ends, e.g., the congestion 
control \cite{ja1988} and real-time monitoring \cite{wi2017}.

\subsection{Dual Problem}

We first reformulate SLRR model \eqref{model3} as
\begin{equation}\label{P}
\begin{aligned}
&\min\limits_{X}\quad \|X\|_{*}+\alpha\|Z\|_F^2\\
&\text{s.t.}\quad  \sum\limits_{j=1}^S R_jXe_j=L,\quad 
\mathcal{P}_{\Omega}(X)=0,\quad X=Y,\quad X-A-Z=0,\quad Y\ge0.
\end{aligned}
\end{equation}
In model \eqref{P}, it follows \begin{equation}\label{aA}
\alpha:=\rho_1+\rho_2\quad \text{and}\quad 
A=\frac{\rho_1\bar{X}+\rho_2\hat{X}}{\rho_1+\rho_2}.
\end{equation}
The Lagrangian function of \eqref{P} is
\begin{multline*}
L(X,Y,Z;U,V,W,Q,G):=\|X\|_{*}+\alpha\|Z\|_F^2-\langle 
U,\mathcal{P}_{\Omega}(X)\rangle-\langle V,X-Y\rangle\\
-\langle W,X-A-Z\rangle-\langle Q,\sum\limits_{j=1}^S R_jXe_j-L\rangle,
\end{multline*}
with the variable matrix $Y$ being nonnegative. 
Since the objective function of \eqref{P} are convex and all of  constraints 
are linear, it has the strong duality, namely
\begin{multline*}
\min\limits_{X,Y\ge0,Z}\max\limits_{U,V,W,G}L(X,Y,Z;U,V,W,Q,G)\\
=\max\limits_{U,V,W,G}\min\limits_{X,Y\ge0,Z}L(X,Y,Z;U,V,W,Q,G).
\end{multline*}
This property leads to the dual problem that is very easy to solve with the 
SCB-spADMM. Next, we derive the dual form through the following minimization 
problem:
$$\min\limits_{X,Y\ge0,Z}L(X,Y,Z;U,V,W,Q,G):=\min\limits_{X,Y\ge0,Z}L_1(X,Y,Z).$$
Note that $X,Y$ and $Z$ in $L_1(\cdot)$ are decoupled. Thus, the above 
problem can be divided into three smaller pieces associated to $X,\ Y$ and $Z$ 
and be solved successively. 

The subproblem of $X$ is solved by 
\begin{multline*}
\min\limits_X\quad  \|X\|_{*}-\langle 
\mathcal{P}_{\Omega}(U)+V+W+\sum\limits_{j=1}^s R_{j}^{T} Q e_{j}^{T},X 
\rangle\\
=\left\{\begin{aligned}
0&,\quad\|\mathcal{P}_{\Omega}(U)+V+W+\sum\limits_{j=1}^S R_{j}^{T} Q 
e_{j}^{T}\|_2\le1,\\
-\infty&,\quad \text{otherwise}.
\end{aligned}\right.
\end{multline*}
The last equality holds according to the inequality $\langle A,B\rangle\le 
\|A\|_{*}\|B\|_2.$ When nonzero $B$ is fixed, the equality can be   attained 
for certain $A$ because $\|\cdot\|_{*}$ and $\|\cdot\|_2$ are dual norms.

In terms of $Y$, the subproblem can be easily fixed by
\begin{equation*}
\begin{aligned}
\min\limits_{Y\ge0}\quad & \langle V,Y\rangle= \left\{\begin{aligned}
0&,\quad V\ge 0,\\
-\infty&,\quad \text{otherwise}.
\end{aligned}\right.
\end{aligned}
\end{equation*}

The optimum of $Z$ can also be trivially obtained since we only need to solve a 
convex quadratic programming 
\begin{equation*}
\begin{aligned}
\min\limits_{Z}\quad & \alpha\|Z\|_F^2+\langle W,Z\rangle= 
-\frac{\|W\|_F^2}{4\alpha}.
\end{aligned}
\end{equation*}

As a consequence, the dual model of \eqref{P} takes the following form, and 
with the strong duality, they are equivalent.
\begin{equation}\label{DD}\begin{aligned}
&\min\limits_{U, V, W, Q, G}\quad \frac{1}{4\alpha}\left\|W-2\alpha 
A\right\|_{F}^{2}-\langle Q, L\rangle \\
&\text { s.t. }\quad \mathcal{P}_{\Omega}(U)+V+W+\sum\limits_{j=1}^S R_{j}^{T} 
Q e_{j}^{T}=G,\quad V\ge0,\quad \|G\|_2\le1.
\end{aligned}\end{equation}

For later developments, we further rewrite problem \eqref{DD} as
\begin{equation}\label{D}\begin{aligned}
&\min\limits_{U, V, W, Q, G}\quad \delta_{\mathbb{R}^{s\times 
s}_{+}}(V)-\langle Q, 
L\rangle+\delta_{\mathcal{B}}(G)+\frac{1}{4\alpha}\left\|W-2\alpha 
A\right\|_{F}^{2} \\
&\text { s.t. }\quad \mathcal{P}_{\Omega}(U)+V+W+\sum\limits_{j=1}^S R_{j}^{T} 
Q e_{j}^{T}=G,
\end{aligned}\end{equation}
where $\mathcal{B}:=\{X\in\mathbb{R}^{S\times S}:\ \|X\|_2\le1\}$, 
$\mathbb{R}^{s\times s}_{+}:=\{X\in\mathbb{R}^{S\times S}:\ X\ge0\}$, and  
$\delta_\mathcal{C}(\cdot)$ is the indicator function on a closed convex set 
$\mathcal{C}$, i.e., $\delta_\mathcal{C}(x)=0$ if $x\in \mathcal{C}$, otherwise 
$\delta_\mathcal{C}(x)=+\infty$. The 
objective of \eqref{D} contains two blocks, namely 
$$\text{block}_1:=\left(\delta_{\mathbb{R}^{s\times s}_{+}}(V)-\langle Q, 
L\rangle\right)\ \text{and}\ 
\text{block}_2:=\left(\delta_{\mathcal{B}}(G)+\frac{1}{4\alpha}\left\|W-2\alpha 
A\right\|_{F}^{2}\right),$$ each of which is a convex composite quadratic  
function. In general, a convex composite quadratic  
function takes  form as
$$f(x,y)=p(x)+\langle y,My+b\rangle,$$ where 
$p(\cdot):\mathcal{X}\to(-\infty,+\infty]$ is closed proper convex function, 
$M$ is a positive semidefinite matrix and $b$ is a vector.

\subsection{Optimality  Condition and Stopping Criterion of 
Algorithm}

The Karush-Kuhn-Tucker (KKT) conditions associated with  \eqref{P} and 
\eqref{D} are given as follows.

\begin{equation*}
\left\{\begin{aligned}
&\sum\limits_{j=1}^S R_jXe_j=L,\quad \mathcal{P}_{\Omega}(X)=0,\quad 
\mathcal{P}_{\Omega}(U)+V+W+\sum\limits_{j=1}^S R_{j}^{T} Q e_{j}^{T}=G,\\
& V\in\mathbb{R}^{S\times S}_{+},\quad G\in\mathcal{B}.
\end{aligned}\right.
\end{equation*}

We can set the stopping criterion according to these KKT conditions. 
Specifically, define $\eta:=\max\{\eta_{P_1},\ \eta_{P_2},\ \eta_D,\ \eta_V,\ 
\eta_G\},$ where 
\begin{equation*}
\left\{\begin{aligned}
&\eta_{P_1}:=\frac{\|\sum\limits_{j=1}^S R_jXe_j-L\|_F}{1+\|L\|_F},\quad 
\eta_{P_2}:=\frac{\|\mathcal{P}_{\Omega}(X)\|_F}{1+\|X\|_F},\\
&\eta_D:=\frac{\|\mathcal{P}_{\Omega}(U)+V+W+\sum\limits_{j=1}^S R_{j}^{T} Q 
e_{j}^{T}-G\|_F}{1+\|G\|_F},\\
&\eta_V:=\frac{\|V-\mathcal{P}_{\mathbb{R}^{S\times 
S}_{+}}(V)\|_F}{1+\|V\|_F},\quad
\eta_G:=\frac{\|G-\mathcal{P}_{\mathcal{B}}(G)\|_F}{1+\|G\|_F}.
\end{aligned}\right.
\end{equation*}
We can set a small positive parameter $\varepsilon$ and stop the algorithm when 
$\eta<\varepsilon$ or the iteration step attains it maximum.

\section{Schur Complement Based Semi-Proximal ADMM}
\label{algorithm}
Next, we would develop the Schur Complement Based semi-proximal ADMM 
(SCB-spADMM) for \eqref{D}. As it is well known,  the ADMM-type algorithm is 
effective and efficient for large-scale convex modes with linear constraints 
\cite{bo2010,chen2018,chen2017,li2014} and has been successfully applied to 
many nuclear norm minimization problems  \cite{wen2012,recht2010,lu2018}. It is 
easy to implement and only requires little memory space. The convergence of our 
SCB-spADMM can be similarly proved as in \cite{li2014}. In particular, the 
optimum of each variable can be solved explicitly.


\subsection{Traditional Multi-Block ADMM}
A natural idea is to directly employ a 5-block ADMM for solving model 
\eqref{D}. Denote 
$$\Gamma(U,V,W,Q,G):=\mathcal{P}_{\Omega}(U)+V+W+\sum\limits_{j=1}^S R_{j}^{T} 
Q e_{j}^{T}-G.$$ The Lagrangian function of \eqref{D} is 
\begin{multline*}
L_\beta(U,V,W,Q,G;X):= \delta_{\mathbb{R}^{S\times 
        S}_{+}}(V)-\langle Q, 
L\rangle+\delta_{\mathcal{B}}(G)+\frac{1}{4\alpha}\left\|W-2\alpha 
A\right\|_{F}^{2}\\
+\langle 
X,\Gamma(U,V,W,Q,G)\rangle+\frac{\beta}{2}\|\Gamma(U,V,W,Q,G)\|_F^2.
\end{multline*}

A commonly used algorithm for solving \eqref{D} is the directly extended 
5-block ADMM. Its iteration steps  are given by
\begin{equation*}
\left\{\begin{aligned}
&U^{k+1}:=\arg\min_{U}L_{\beta}(U,V^{k},W^{k},Q^{k},G^{k};X^{k}),\\
&V^{k+1}:=\arg\min_{V}L_{\beta}(U^{k+1},V,W^{k},Q^{k},G^{k};X^{k}),\\
&W^{k+1}:=\arg\min_{W}L_{\beta}(U^{k+1},V^{k+1},W,Q^{k},G^{k};X^{k}),\\
&Q^{k+1}:=\arg\min_{Q}L_{\beta}(U^{k+1},V^{k+1},W^{k+1},Q,G^{k};X^{k}),\\
&G^{k+1}:=\arg\min_{G}L_{\beta}(U^{k+1},V^{k+1},W^{k+1},Q^{k+1},G;X^{k}),\\
&X^{k+1}:=X^{k}+\tau\beta(\mathcal{P}_{\Omega}(U^{k+1})+V^{k+1}+W^{k+1}+\sum\limits_{j=1}^s
 R_{j}^{T} Q^{k+1} e_{j}^{T}-G^{k+1}).
\end{aligned}\right.
\end{equation*}
Since each variable can be obtained explicitly in the iterations, the solving 
of subproblems becomes very simple. On the other hand, as mentioned in 
\cite{chen2016}, it is a challenge to ensure that the designed multi-block 
($\ge3$ blocks) ADMM is convergent. This is a difficulty we need to overcome.

Inspired by Li, Sun and Toh's impressive work \cite{li2014}, we will propose a 
convergent ADMM algorithm below. By introducing the proximal or semi-proximal 
terms, the subproblems in the proximal ADMM (pADMM) or semi-proximal ADMM 
(spADMM) can be solved efficiently. Thus they are more suitable to cope with 
convex composite quadratic programmings and convex quadratic conic 
programmings. The key point of designing the semi proximal term in SCB-spADMM  
is to modify the algorithm into the split version based on Gauss-Seidel type 
decomposition.

\begin{algorithm}[h]
    \renewcommand{\thealgorithm}{}
    \caption{SCB-spADMM for solving \eqref{D}}
    \label{SCB-SPADMM}
    \begin{algorithmic}
        \STATE{Define $R_j$ as the $j-$th block of the routing matrix $R$, and 
            $e_j$ as the zero vector except the $j$-th element being one. 
            Calculate 
            $\alpha$ and $A$ through \eqref{aA}. Choose $U^0,\ Q^0,\ V^0,\ 
            W^0,\ 
            G^0,\ X^0$ as zero matrices. For $k=0,1,\cdots,$ repeat the 
            following 
            iteration scheme.}
        \WHILE{Stopping criterion is not reached}
        \STATE{Update $U,\ Q$ and $V$ as the following order according to 
            formulations \eqref{U}, \eqref{Q}, \eqref{V}, respectively.			
            \hspace{0.6in}
            \begin{equation}\label{P1}
            \left\{\begin{aligned}
            &U^{k+\frac{1}{2}}:=\arg\min\limits_UL_{\beta}(U,Q^{k},V^{k},W^{k},G^{k};X^{k})+\frac{\beta}{2}\|U-U^k\|_{\mathcal{H}_U}^2,\\
            &Q^{k+\frac{1}{2}}:=\arg\min\limits_QL_{\beta}(U^{k+\frac{1}{2}},Q,V^{k},W^{k},G^{k};X^{k})+\frac{\beta}{2}\|Q-
            Q^k\|_{\mathcal{H}_Q}^2,\\
            &V^{k+1}:=\arg\min\limits_VL_{\beta}(U^{k+\frac{1}{2}},Q^{k+\frac{1}{2}},V,W^{k},G^{k};X^{k})+\frac{\beta}{2}\|V-
            V^k\|_{\mathcal{H}_V}^2,\\
            &Q^{k+1}:=\arg\min\limits_QL_{\beta}(U^{k+\frac{1}{2}},Q,V^{k+1},W^{k},G^{k};X^{k})+\frac{\beta}{2}\|Q-
            Q^{k+\frac{1}{2}}\|_{\mathcal{H}_Q}^2,\\
            &U^{k+1}:=\arg\min\limits_UL_{\beta}(U,Q^{k+1},V^{k+1},W^{k},G^{k};X^{k})+\frac{\beta}{2}\|U-U^{k+\frac{1}{2}}\|_{\mathcal{H}_U}^2,\\
            \end{aligned}\right.
            \end{equation}			
            Update $W$ and $G$ as the following order according to formulations 
            \eqref{W}, \eqref{G}, respectively.			
            \hspace{0.6in}
            \begin{equation}\label{P2}
            \left\{\begin{aligned}
            &W^{k+\frac{1}{2}}:=\arg\min\limits_WL_{\beta}(U^{k+1},Q^{k+1},V^{k+1},W,G^{k};X^{k})+\frac{\beta}{2}\|W-
            W^k\|_{\mathcal{H}_W}^2,\\
            &G^{k+1}:=\arg\min\limits_GL_{\beta}(U^{k+1},Q^{k+1},V^{k+1},W^{k+\frac{1}{2}},G;X^{k})+\frac{\beta}{2}\|G-
            G^k\|_{\mathcal{H}_G}^2,\\
            &W^{k+1}:=\arg\min\limits_WL_{\beta}(U^{k+1},Q^{k+1},V^{k+1},W,G^{k+1};X^{k})+\frac{\beta}{2}\|W-
            W^{k+\frac{1}{2}}\|_{\mathcal{H}_W}^2,\\
            \end{aligned}\right.
            \end{equation}			
            $X^{k+1}:=X^{k}+\gamma\beta(P_{\Omega}(U^{k+1})+V^{k+1}+W^{k+1}+\sum\limits_{j=1}^s
            R_{j}^{T} Q^{k+1} e_{j}^{T}-G^{k+1}).$			
        }
        \ENDWHILE
        \RETURN $X$
    \end{algorithmic}
\end{algorithm}

\subsection{Algorithm Framework}
As introduced in \cite{li2014}, 
let $\mathcal{H}_U,\ \mathcal{H}_Q,\ \mathcal{H}_W,\ \mathcal{H}_V$ and  
$\mathcal{H}_G$ be self-adjoint positive semidefinite linear operators and 
satisfy
\begin{equation}\label{H}
\begin{aligned}
&\mathcal{H}_U:=\mathcal{G}_U-P_{\Omega}^{*}P_{\Omega}=\mathcal{G}_U-P_{\Omega},\quad
\mathcal{H}_Q:=\mathcal{G}_Q-\mathcal{A}^{*}\mathcal{A},\\
&\mathcal{H}_W:=\mathcal{G}_W-\frac{1}{\beta}\mathcal{I}-\mathcal{I}^{*}\mathcal{I}=\mathcal{G}_W-(\frac{1}{\beta}+1)\mathcal{I},\\
&\mathcal{H}_V:=\mathcal{G}_V-\mathcal{I}^{*}\mathcal{I}=\mathcal{G}_V-\mathcal{I},\quad\mathcal{H}_G:=\mathcal{G}_G-\mathcal{I}^{*}\mathcal{I}=\mathcal{G}_G-\mathcal{I},\\
\end{aligned}
\end{equation}
where $\mathcal{G}_U,\ \mathcal{G}_Q,\ \mathcal{G}_W,\ \mathcal{G}_V$ and 
$\mathcal{G}_G$ are also self-adjoint positive semidefite linear operators, and 
additionally, their inverses are relatively easy to compute. 
$\mathcal{A}$ is a linear operator satisfying 
$\mathcal{A}Q:=\sum\limits_{j=1}^s R_{j}^{T} Q e_{j}^{T}$ and $\mathcal{I}$ is 
the identity operator, i.e., $\mathcal{I}X=X$ for any $X$. From the definition, 
we can derive that $\mathcal{A}^{*}\mathcal{A}Q=\sum\limits_{j=1}^sR_jR_j^TQ.$
We propose the SCB-spADMM algorithm for solving model \eqref{D} as in Algorithm 
\ref{SCB-SPADMM}.

\subsection{Convergence}
The presence of $\mathcal{H}_U,\ \mathcal{H}_Q,\ \mathcal{H}_W,\ \mathcal{H}_V$ 
and  $\mathcal{H}_G$ can guarantee the existence of the solution of subproblem 
\eqref{P1} and \eqref{P2}, as well as the boundedness of the iteration points 
$U^k,\ Q^k,\  W^k,\ V^k$ and $G^{k}$. In general, $\mathcal{H}_U,\ 
\mathcal{H}_Q,\ \mathcal{H}_W,\ \mathcal{H}_V$ and  $\mathcal{H}_G$ should be 
as small as possible  and meanwhile $U^k,\ Q^k,\  W^k,\ V^k$ and $G^{k}$ are 
still easy to compute. According to \cite{li2014}, we have the following 
convergence result.

\begin{Theorem}
    If one of (1) $\tau\in(0,\frac{\sqrt{5}+1}{2})$ and (2) 
    $\tau\ge\frac{\sqrt{5}+1}{2}$ but 
    $\sum\limits_{k=0}^{+\infty}\|(G^{k+1}-G^{k})+(W^{k+1}-W^{k})\|_F^2+\frac{1}{\tau}\|\mathcal{P}_{\Omega}(U^{k+1})+V^{k+1}+W^{k+1}+\sum\limits_{j=1}^s
     R_{j}^{T} Q^{k+1} e_{j}^{T}-G^{k+1}\|_F^2<+\infty$, then we have the 
    following properties.
    \begin{enumerate}[(i)]
        \item  $U^k,$ $Q^k$,  $W^k$, $V^k$ and $G^{k}$ are automatically well 
        defined based on the settings of $\mathcal{H}_U$, $\mathcal{H}_Q,\ 
        \mathcal{H}_W,$ $\mathcal{H}_V$ and  $\mathcal{H}_G$ given by \eqref{H}.
        \item $\{(U^{k},\ Q^{k},\ W^{k}$, $V^{k},\ G^{k},\ 
        X^k)\}_{k=1}^{\infty}$ converges to the unique accumulate point 
        $(U^{\infty},\ Q^{\infty},\ W^{\infty},\ V^{\infty},$ $ G^{\infty},\ 
        X^{\infty})$, with $(U^{\infty},\ Q^{\infty},\ W^{\infty},\ 
        V^{\infty},\ G^{\infty})$ being the optimum of model \eqref{D} and 
        $X^{\infty}$ being the optimum of model \eqref{P}.
    \end{enumerate}
\end{Theorem}

\begin{proof}
    After substituting the specific form of \eqref{D}, it becomes a 
    straightforward deduction of Theorem 3 in \cite{li2014}. $\Box$
\end{proof}

Recall the requirements of $\mathcal{H}_U,\ \mathcal{H}_Q,\ \mathcal{H}_W,\ 
\mathcal{H}_V$ and  $\mathcal{H}_G$ introduced previously, we can simply set  
$\mathcal{H}_W=\mathcal{H}_V=\mathcal{H}_G={\bf 0}$ if we let 
$\mathcal{G}_V=\mathcal{G}_G=(\frac{1}{\beta}+1)^{-1}\mathcal{G}_W=\mathcal{I}.$
Since $\mathcal{P}_{\Omega}$ is a nondegenerate projection operator, it 
follows that $\|\mathcal{P}_{\Omega}\|_2=1$. Thus, we can set 
$\mathcal{G}_U=\mathcal{I}$ and have 
$\mathcal{H}_U=\mathcal{I}-\mathcal{P}_{\Omega}.$ In terms of $\mathcal{H}_Q$, 
from the equation $\mathcal{A}^{*}\mathcal{A}Q=\sum\limits_{j=1}^SR_jR_j^TQ$, 
we know that  
$\|\mathcal{A}^{*}\mathcal{A}\|_2=\|\sum\limits_{j=1}^SR_jR_j^T\|_2.$ As for 
$\mathcal{G}_Q$, a common and simple setting is  
$\mathcal{G}_Q=\lambda_{\text{max}}\mathcal{I}$, where $\lambda_{\text{max}}$ 
is the largest eigenvalue of $\mathcal{A}^{*}\mathcal{A}$. Thus, we have that 
$$\mathcal{H}_Q(Q)=(\lambda_{\text{max}}I-\sum\limits_{j=1}^SR_jR_j^T)Q:=H_QQ.$$

In the following, we will derive the explicit forms of the subproblems 
associated to $U,\ Q,\ V,\ W$ and $G$, respectively.

The subproblem of $U$ takes the form as 
$$\min\limits_{U}\ \langle 
X,P_\Omega(U)\rangle+\frac{\beta}{2}\|\Gamma(U,V,W,Q,G)\|_F^2+\frac{\beta}{2}\|U-U^k\|_{\mathcal{I}-\mathcal{P}_\Omega}^2.$$
Since it is an unconstraint model and its objective function is convex and 
differentiable, we can obtain its optimum by directly setting the gradient of 
the objective function to zero. Thus, the optimum of $U$ is
\begin{equation}\label{U}
U^{*}:=\mathcal{P}_{\Omega}\left(-V-\sum\limits_{j=1}^S R_{j}^{T} Q 
e_{j}^{T}-W+G-\frac{X}{\beta}\right)+\mathcal{P}_{\Omega^\text{C}}(U^k),
\end{equation}
Similarly, we can simply derive the optimum of $Q$ and $W$ by
\begin{equation}\label{Q}
\begin{aligned}
Q^{*}:=&\arg\min\limits_Q\ \langle \sum\limits_{j=1}^S 
R_jXe_j-L,Q\rangle+\frac{\beta}{2}\|\Gamma(U,V,W,Q,G)\|_F^2+\frac{\beta}{2}\|Q-Q^k\|_{H_Q}^2\\
=&-\frac{1}{\lambda_\text{max}}\left(\sum\limits_{j=1}^S 
R_{j}(V+P_{\Omega}(U)+W-G) e_{j}-H_QQ^k+\frac{\sum\limits_{j=1}^s R_{j} X 
    e_{j}-L}{\beta}\right)
\end{aligned}
\end{equation}
and
\begin{equation}\label{W}
\begin{aligned}
W^{*}:=&\arg\min\limits_W\ \frac{1}{4\alpha}\|W-2\alpha A\|_F^2+\langle 
X,W\rangle+\frac{\beta}{2}\|\Gamma(U,V,W,Q,G)\|_F^2\\
=&\frac{A-X-\beta\left(V+\sum\limits_{j=1}^S R_{j}^{T} Q 
e_{j}^{T}+P_{\Omega}(U)-G\right)}{\frac{1}{2 \alpha}+\beta}.
\end{aligned}
\end{equation}
The subproblem in terms of $V$ can be divided into a group of convex quadratic 
programmings, each of which asks for nonnegative solution and associates to 
scalar variable. It is easy to derive that
\begin{equation}\label{V}
\begin{aligned}
V^{*}:=&\arg\min\limits_V\ \delta_{\mathbb{R}^{S\times S}_{+}}(V)+\langle 
X,V\rangle+\frac{\beta}{2}\|\mathcal{P}_{\Omega}(U)+V+W+\sum\limits_{j=1}^s 
R_{j}^{T} Q e_{j}^{T}-G\|_F^2\\
=&\max \left\{0,-\left(\sum_{j=1}^s R_{j}^{T} Q 
e_{j}^{T}+P_{\Omega}(U)+W-G\right)-\frac{X}{\beta}\right\}.
\end{aligned}
\end{equation}
About $G$, suppose $\ U_{G} \Sigma_{G} V_{G}:=V+\sum R_{j}^{T} Q 
e_{j}^{T}+P_{\Omega}(U)+W-\frac{X}{\beta}$ is the economical singular value 
decomposition ({\sl SVD}). Then, we have that
\begin{equation}\label{G}
\begin{aligned}
G^{*}:=&\arg\min\limits_G\ \delta_{\mathcal{B}}(G)-\langle 
X,G\rangle+\frac{\beta}{2}\|\Gamma(U,V,W,Q,G)\|_F^2\\
=&U_{G} \min \left\{I, \Sigma_{G}\right\} V_{G}^{T}.
\end{aligned}
\end{equation}

\section{Numerical Simulations}\label{ns}

In this section, we present the simulation results of our SLLR model and 
SCH-spADMM algorithm on the Abilene and GEANT datasets. The results of 
Tomo-SRMF \cite{ro2012}, SDPLR \cite{recht2010} and Tomo-Gravity 
\cite{zhang2005}  
are also output for  
comparison. As already mentioned in the introduction section, we will focus on 
more practical network tomography problem 
\cite{TMLink04,recht2010,ro2012,tmverdi96,tmcao00,TMLink03} here rather than 
the completion  problem 
\cite{ca2010,ch2012,ji2020,recht2010,ro2012,xie2017,2019Active,xie2018,xie2015,xie2019,zh2015}.

\begin{figure}[!h]
    \centerline{\includegraphics[width=16cm,height=3.8cm]{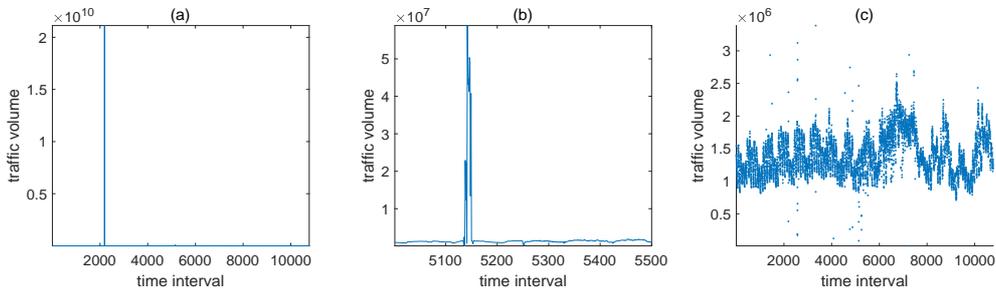}}
    \caption{The Frobenius norm of (a) the original G$\rm\E$ANT traffics in 
    each 
        time interval,
        (b) the original G$\rm\E$ANT traffics from the 5001st 
        time interval to the 5500th time interval, (c) the 
        modified G$\rm\E$ANT traffics in each time interval.}\label{ga}
\end{figure}
\begin{figure}[!h]
    \centerline{\includegraphics[width=15cm,height=6cm]{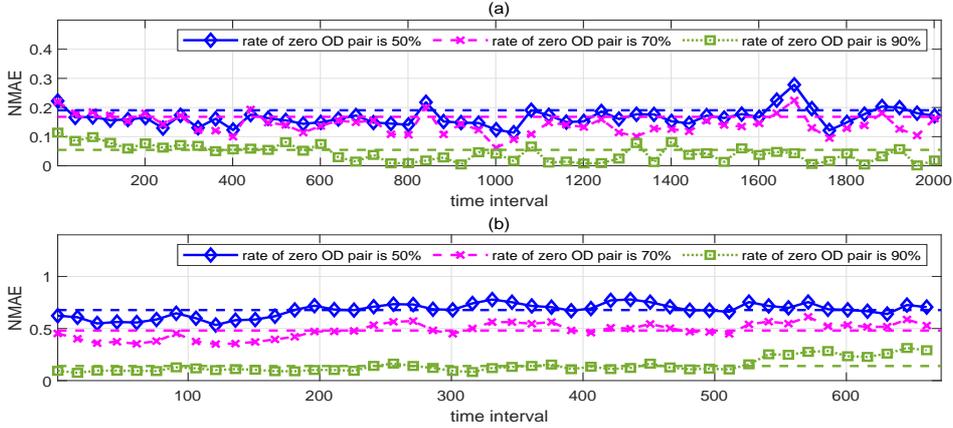}}
    \caption{NMAE in different time intervals on (a) the complete three-week 
        traffics of Abilene dataset, (b) the first-week traffics of
        G$\rm\E$ANT dataset. The spaces between adjacent time intervals are 40 
        for 
        Abilene dataset and 15 for G$\rm\E$ANT dataset.}\label{timeslot}
\end{figure}

We notice there are some anomalies in G$\rm\E$ANT traffics \cite{uh2006}. 
As 
seen in Figure \ref{ga}, a few of traffic volumes are much larger than the 
others. 
Suppose the traffics in the time intervals $t_1$ and $t_2$ are correct but all 
the intermediate traffics are abnormal. Then, we modify 
the anomalies by using linear interpolation between the time intervals $t_1$ 
and $t_2$. In Figure \ref{ga} (c), the refined data is presented. For Abilene 
dataset, we can see that it has good quality from Figure \ref{timeslot}. Thus 
no special processing is required.

As previously mentioned, the real-world traffic data is of high sparsity.
In order to imitate this realistic circumstance, we modify the Abilene dataset 
and G$\rm\E$ANT dataset by setting the smallest $p\%\ (p=50,70,90)$ of OD pairs 
to zero-ODs, which means they do not transmit any network traffic.

K-fold \cite{mc2004} Cross Validation (CV) is employed to find proper 
parameters 
$(\alpha_1, \alpha_2, \beta)$ in our algorithm. We uniformly divide the rows 
of the link-measurement matrix $L$ into $K$ groups. Each row group and the 
traffics therein are denoted as $\Psi_l$ and $L_{\Psi_l}$ $(l=1,\cdots,K)$,  
respectively. Then, we generate candidate parameters $\{(\alpha_1^j, 
\alpha_2^j, \beta^j)\}_{j=1}^{1000}$ from a wide-range interval. 
All of candidates will be evaluated by the following steps:
\begin{enumerate}
    \item Let the integer $l$ change from 1 to $K$. In each case, $L_{\Psi_l}$ 
    is taken as the test set, while the other link-load data is taken as the 
    training set.
    \item Based on the link-load traffics in the training set, we can obtain an 
    estimation of the traffics in rows $\Psi_l$ using proposed SCB-spADMM, 
    denoted as $\hat{L}_{\Psi_l}$.
    \item If $l=K$, calculate the overall recovery error by Normalized 
    Mean Absolute Error ({NMAE}) for CV, defined as 
    \begin{equation*}\label{NMAE1}
    \text{N}_{\text{CV}}:=\frac{\sum\limits_{l=1}^K\sum\limits_{i\in\Psi_l}\sum\limits_{j=1}^T\left|(\hat{L}_{\Psi_l})_{ij}-(L_{\Psi_l})_{ij}\right|}{\sum\limits_{i=1}^M\sum\limits_{j=1}^TL_{ij}}.
    \end{equation*}
\end{enumerate}
Ultimately, we choose the candidate that corresponds to the lowest 
$\text{N}_{\text{CV}}$ as the parameters used in our algorithm. 
With the fixed parameters, we can implement SCB-spADMM to estimate the whole OD 
traffics. The recovery performance is evaluated by NMAE on the missing data:
\begin{equation}\label{SNMAE}
\text{NMAE}:=\frac{\sum\limits_{(i,j,k)\in\Omega^C}|\hat{\mathcal{X}}_{ijk}-\mathcal{X}_{ijk}|}{\sum\limits_{(i,j,k)\in\Omega^C}\mathcal{X}_{ijk}},
\end{equation}

We compare SLRR against three alternative methods: Tomo-SRMF, SDPLR 
and Tomo-Gravity. The complete three-week Abilene traffics and four-month 
G$\rm\E$ANT traffics are simulated.  In particular, we simulate G$\rm\E$ANT 
dataset in three different time periods: the first week, the first eight weeks 
and the last eight weeks, respectively.  We need to point out that other work  
\cite{ro2012,zh2015} only simulated the one-week data of the GEANT dataset. In 
fact, sometimes people only care about one-week traffics rather than the 
aggregate data.

The simulation results are reported in Table \ref{AbileneTest}. We can see that 
as the sparsity increases, the recovery errors will decrease. This is because 
high sparsity contains more raw traffic information and reduce uncertainty. 
More importantly, our SLRR is significantly better than other methods in all 
cases. For the GEANT dataset, the results of the first-week data are better 
than the results of the other two groups. This may be related to the data 
quality of the GEANT dataset. Nevertheless, our SLRR is still superior to other 
methods in these cases.

\begin{table}[!h]
    \centering
    \caption{Comparison of NMAE on Abilene and G$\rm\E$ANT 
        datasets. For G$\rm\E$ANT dataset, we output the simulation results of 
        the 
        first-week/the first-eight-week/the last-eight-week data.}  
    \label{AbileneTest}
    \begin{tabular}{|c|c|c|c|}  
        \hline
        Method & Sparsity = 50\% & Sparsity = 70\% & Sparsity = 90\%  \\  
        \hline
        & \multicolumn{3}{c|}{Abilene}\\
        \hline
        Tomo-SRMF & 0.374 & 0.335 & 0.236\\ 
        SDPLR & 0.535 & 0.393 & 0.272 \\
        Tomo-Gravity & 0.631 & 0.558 & 0.248\\ 
        SLRR  & {\bf 0.193} & {\bf 0.136} & {\bf 0.047}\\
        \hline
        & \multicolumn{3}{c|}{G$\rm\E$ANT}\\
        \hline
        Tomo-SRMF & 1.003/1.182/1.342 & 0.725/0.922/1.075 & 
        0.271/0.465/0.530\\
        SDPLR & 1.004/1.055/1.171 & 0.734/0.859/0.973 & 
        0.308/0.456/0.526\\
        Tomo-Gravity& 1.065/1.113/1.145 & 0.871/0.971/1.005 & 
        0.583/0.644/0.734\\
        SLRR  & {\bf 0.640/0.743/0.792} & {\bf 0.444/0.590/0.639} & {\bf 
            0.124/0.271/0.337}\\
        \hline
    \end{tabular}  
\end{table}

The recovery performances of SLRR in different time intervals are presented in 
Figure \ref{timeslot}. We especially mark the NAME under certain sparsity with 
a 
horizontal dashed line of the same color. It shows that our recovery results 
are consistent and reliable.

\section{The HOD Dataset and Simulations}\label{hod}

\begin{figure}[!h]
    \centerline{\includegraphics[width=15cm,height=11cm]{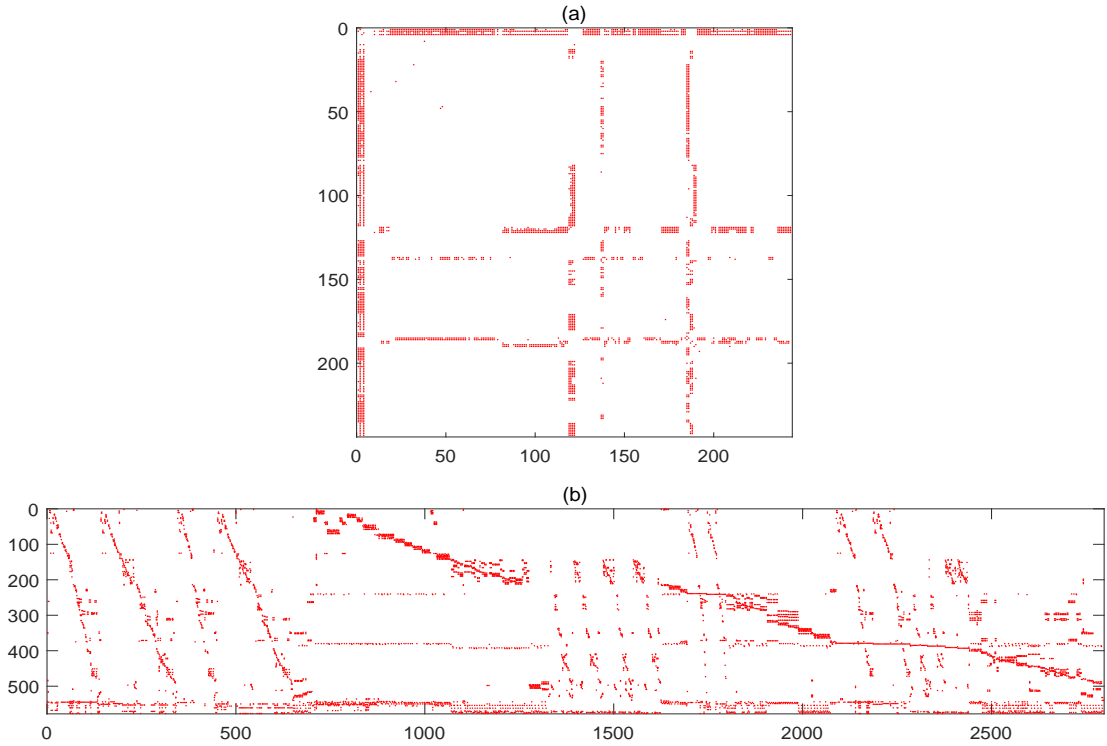}}
    \caption{The sparsity of (a) the spatial matrix, (b) the routing 
        matrix 
        of HOD dataset. The red components in the spatial matrix means that the 
        corresponding OD pairs have transmitted traffic data, and those in 
        the routing matrix stand for 1-components.}\label{HOD}
\end{figure}

\begin{figure}[!h]
    \centerline{\includegraphics[width=15cm,height=5.5cm]{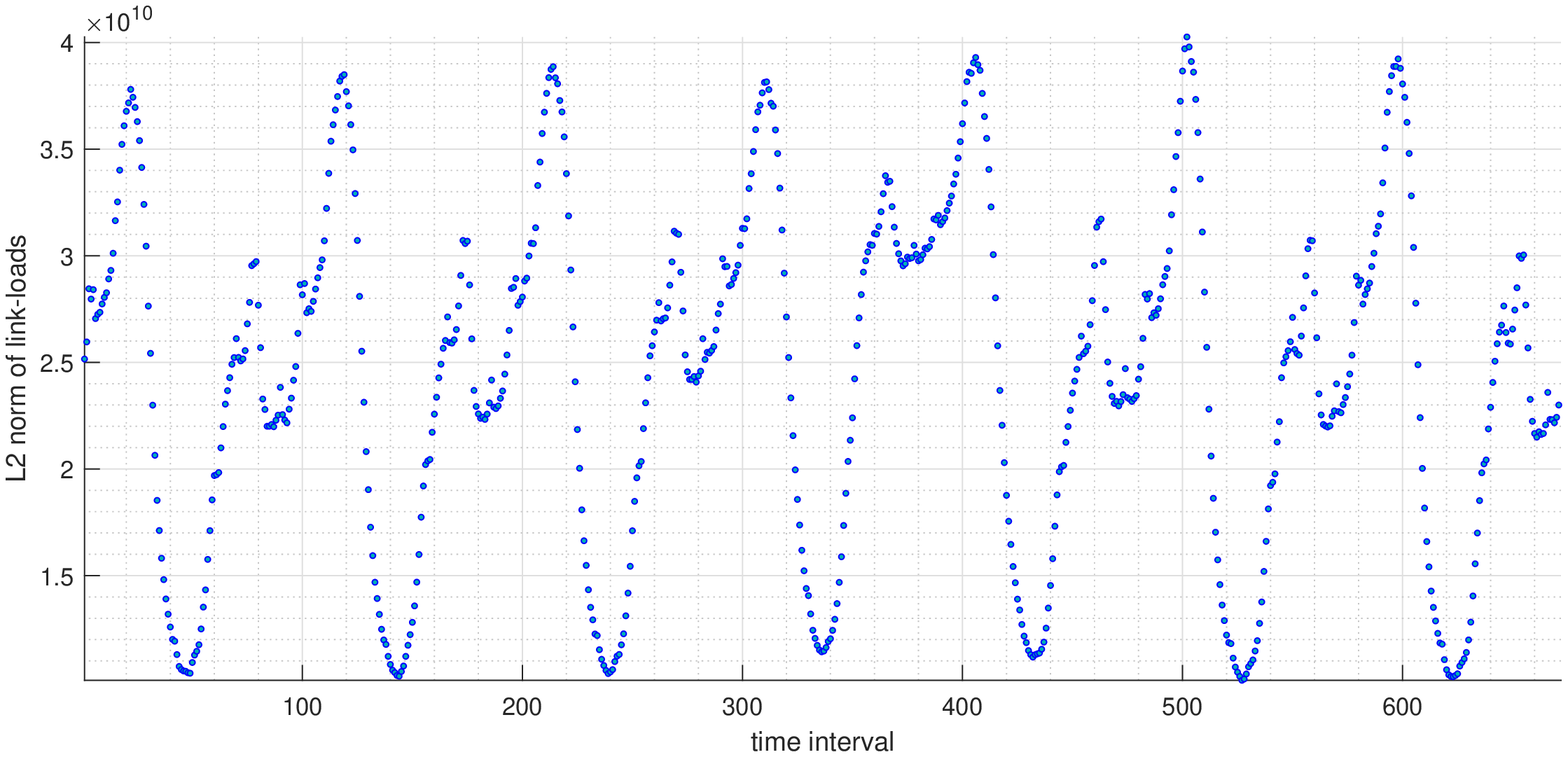}}
    \caption{The Frobenius norm of link-loads in different time intervals. The 
        figure illustrates the the continuity and the periodicity of HOD 
        dataset.}\label{linkload}
\end{figure}

Compared with Abilene dataset and GEANT dataset, the most significant feature 
of HOD dataset is that it is highly sparsity. In fact, 95.26\% of OD pairs have 
no traffic. In Figure \ref{HOD}, we fully demonstrate this. Moreover, HOD 
dataset 
also has continuity and periodicity, see Figure \ref{linkload} for 
illustration. 
In the following simulations, we will also make use of these properties.

Next, we evaluate different methods on HOD dataset. Because the real OD 
traffics are unclear, NMAE metric \eqref{SNMAE} is no longer effective. 
Again we utilize CV errors \eqref{NMAE1} for assessment. Except the 
aforementioned K-fold CV, we also consider another kind of CV, named the Monte 
Carlo CV \cite{du2007} to make fair evaluation. In the Monte 
Carlo CV, the training set and the test set are 
randomly generated. We set the ratio of the test sets to 2\% and repeat 50 
times in each simulation.

Table \ref{IndTest} presents the recovery performance of different methods 
testing on HOD dataset. `SLRR-MC' stands for the case where SLRR is evaluated 
by 
Monte Carlo CV.
From Table \ref{IndTest}, we find that over the four methods, SLRR  
significantly outperforms the others. Moreover, the recovery results of 
SLRR assessed by two different CVs are very close. 
\begin{table}[!h]
    \centering
    \caption{Comparison of recovery performance on HOD dataset.}  
    \label{IndTest}
    \begin{tabular}{|c|c|c|c|c|c|}  
        \hline
        Method & Tomo-SRMF & SDPLR & Tomo-Gravity & SLRR & SLRR-MC  \\ 
        \hline 
        NMAE & 0.3314 & 0.2832 & 0.3243 & {\bf 0.1588} & {\bf 0.1519} \\
        \hline
    \end{tabular}  
\end{table}

\begin{figure}[!h]
    \centerline{\includegraphics[width=15cm,height=7.5cm]{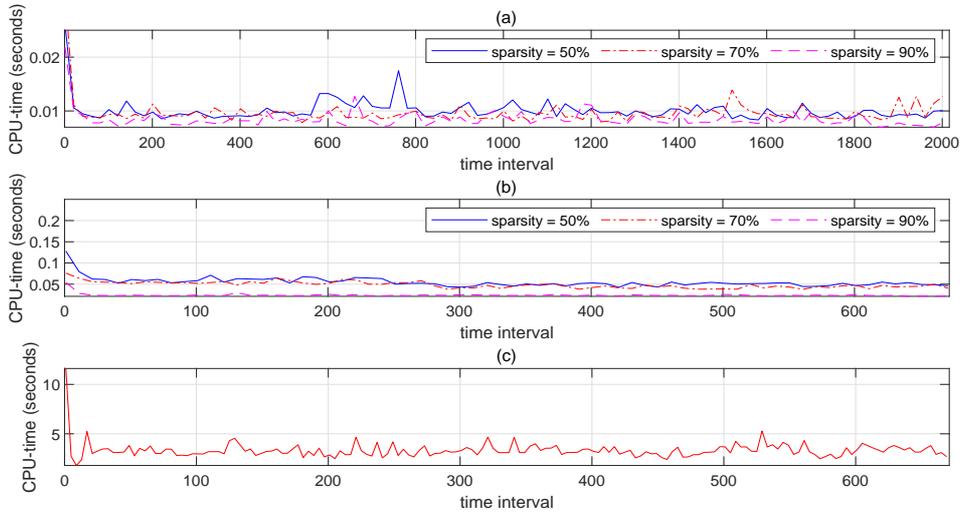}}
    \caption{CPU-time of each time interval on (a) Abilene dataset, (b) 
        G$\rm\E$ANT dataset, (c) HOD dataset. For Abilene and G$\rm\E$ANT 
        datasets, 
        it takes less than 0.1 second to recover the traffics in one time 
        interval. 
        For large-scale HOD dataset, it also only takes several seconds for 
        each 
        interval.}\label{cputime}
\end{figure}

The computational time of our algorithm is also investigated. The overall 
simulations are implemented on an ASUS laptop (4-core, Intel i7-10510U, 
@2.30GHz, 16G RAM). 
We have to admit that our SLRR model is slower than the other methods due to 
the SVD decomposition. However, our method only takes a few seconds (see 
Figure  
\ref{cputime}) to recover the time interval of more than ten minutes in 
practical scenarios. This well satisfies the requirement of practical problems. 
Moreover, our SLRR model has better recovery accuracy than all the existing 
methods in any case. These make our SLRR model have significant advantages in 
the network tomography problem, and it may be integrated into future products.

\section{Conclusions}
\label{conc}
Our research extends the knowledge  into solving large-scale network 
tomography problem. By leveraging on the spatial low-rank property and the 
sparsity of traffic data, we propose SLRR 
model to accurately inference the network traffics. The separability of SLRR 
can significantly reduce the problem scale and make it possible to design 
parallel algorithm. To solve SLRR, we design a globally convergent Schur 
complement based semi-proximal alternating direction method of multipliers. 
Since all of optimums of block variables can be calculated by 
explicit forms, the iteration scheme is fast and accurate. K-fold cross 
validation is employed to find proper parameters in our algorithm. By taking 
advantages of the above techniques, our method outperforms all the others with 
respect to the recovery accuracy. Moreover, the seconds-feedback on different 
datasets also demonstrate the high computation efficiency of our algorithm. 
Except the proposed methodology, another highlight of our work is that we 
publish an industrial dataset HOD to enrich the database of network traffics 
recovery study.


\section*{acknowledgement}
    This work was supported by the National Natural Science Foundation of China 
    (Grant No. 11771244,11871297) and Tsinghua University Initiative Scientific 
    Research Program.




\end{document}